\documentclass{article}
\usepackage{amssymb,amsmath,bm,hyperref}
\usepackage{amssymb,amsmath,bm,geometry,graphics,graphicx,url,color}
\def\no{\noindent}
\def\pmatrix{\left(\begin{array}}
\def\endpmatrix{\end{array}\right)}

\def\RR{\mathbb{R}}

\def\B{{\cal B}}
\def\C{{\cal C}}

\def\I{{\cal I}}

\def\P{{\cal P}}

\def\dd{\mathrm{d}}

\def\ii{\mathrm{i}}

\newtheorem{theo}{Theorem}
\newtheorem{lem}{Lemma}
\newtheorem{cor}{Corollary}
\newtheorem{rem}{Remark}

\def\proof{\noindent\underline{Proof}\quad}
\def\QED{\mbox{~$\Box{~}$}}

\def\bfb{{\bm{b}}}
\def\bfc{{\bm{c}}}

\def\bfe{{\bm{e}}}

\def\bfgamma{{\bm{\gamma}}}
\def\bfdelta{{\bm{\delta}}}
\def\bfeta{{\bm{\eta}}}

\def\eps{\varepsilon}

\def\hrho{{\bar\rho}}
\def\hM{{\bar M}}

\def\hg{{\hat\gamma}}
\def\hbfg{{\hat\bfgamma}}

\begin{document}

\title{Analysis of Spectral Hamiltonian Boundary Value Methods (SHBVMs) for the numerical solution of ODE problems}

\author{Pierluigi Amodio\footnote{Dipartimento di Matematica, Universit\`a di Bari, Italy. \url{{pierluigi.amodio,felice.iavernaro}@uniba.it}} \and Luigi Brugnano\footnote{Dipartimento di Matematica e Informatica ``U.\,Dini'', Universit\`a di Firenze, Italy. \url{luigi.brugnano@unifi.it}}  \and  Felice Iavernaro$^*$}


\maketitle

\begin{abstract} Recently, the numerical solution of stiffly/highly-oscillatory Hamiltonian problems has been attacked by using Hamiltonian Boundary Value Methods (HBVMs) as spectral methods in time. While a theoretical analysis of this spectral approach has been only partially addressed, there is enough numerical evidence that it turns out to be very effective even when applied to a wider range of problems. Here we fill this gap by providing a thorough convergence analysis of the methods and confirm the theoretical results with the aid of  a few numerical tests.

\medskip
\no{\bf Keywords:~}Spectral methods, Legendre polynomials, Hamiltonian Boundary Value Methods, HBVMs, SHBVMs.

\medskip
\no{\bf MSC:~} 65L05, 65P10.

\end{abstract}

\section{Introduction}\label{intro} In recent years, the efficient numerical solution of Hamiltonian problems has been tackled via the definition of the energy-conserving Runge-Kutta methods named Hamiltonian Boundary Value Methods (HBVMs) \cite{BIT2009,BIS2010,BIT2010,BIT2011,BIT2012,BIT2012-1,BIT2015}. Such methods have been developed along several directions (see, e.g., \cite{BS2014,BCMR2012,BIT2012-2}), including Hamiltonian BVPs \cite{ABF2015} and Hamiltonian PDEs \cite{BFCI2015,BBFCI2018,BGS2018,BZL2018} (we also refer to the monograph \cite{LIMbook2016} and to the recent review paper \cite{BI2018} for further details). 

More recently, HBVMs have been used as spectral methods in time for solving highly-oscillatory Hamiltonian problems \cite{BMR2018}, as well as stiffly-oscillatory Hamiltonian problems \cite{BIMR2018} emerging from the space semi-discretization of Hamiltonian PDEs. Their spectral implementation is justified by the fact that this family of methods perform a projection of the vector field onto a finite dimensional subspace via a least square approach based on the use of Legendre orthogonal polynomials \cite{BIT2012}. This spectral approach, supported by a very efficient nonlinear iteration technique to handle the large nonlinear systems needed to advance the solution in time  (see  \cite{BIT2011}, \cite[Chapter\,4]{LIMbook2016} and \cite{BI2018}), proved to be very effective.  However, a thorough convergence analysis of HBVMs, used as spectral methods, was still lacking.  In fact, when using large stepsizes, as is required by the spectral startegy, the notion of classical order of a method is not sufficient to explain the correct asymptotic behaviour of the solutions, so that a different analysis is needed, which is the main goal of the present paper. Moreover, the theoretical achievements will be numerically confirmed by applying the methods to a number of  ODE-IVPs. 

It is worth mentioning that early references on the use of spectral methods in time are, e.g., \cite{Hu1972,Hu1972-1,Bo1997,BeSt2000}, while a related reference is \cite{TS2012}. 

With this premises, the structure of the paper is as follows: in Section~\ref{spectime} we analyze the use of spectral methods in time; in Section~\ref{shbvms} we discuss the efficient implementation of the fully discrete method; in Section~\ref{numtest} we provide numerical evidence of the effectiveness of such an approach, confirming the theoretical achievements. At last, a few conclusions are reported in Section~\ref{fine}.

\section{Spectral approximation in time}\label{spectime}

This section contains the main theoretical results regarding the spectral methods that we shall use for the numerical solution of the ODE-IVP
\begin{equation}\label{ivp}
\dot y(t) = f(y(t)), \qquad y(0) = y_0\in\RR^m. 
\end{equation}
Hereafter, $f$ is assumed to be suitably smooth (in particular, we shall assume $f(z)$ to be analytic in a closed complex ball centered at $y_0$). We consider the solution of problem (\ref{ivp}) on the interval $[0,h]$, where $h$ stands for the time-step to be used by a one-step numerical method. The same arguments will be then repeated for the subsequent integration steps.  According to \cite{BIT2012}, we consider the expansion of the right-hand side of (\ref{ivp}) along the shifted and scaled Legendre polynomial orthonormal basis $\{P_j\}_{j\ge0}$,
$$
P_j\in\Pi_j, \qquad \int_0^1 P_i(x)P_j(x)\dd x=\delta_{ij},\qquad i,j=0,1,\dots,
$$
with $\Pi_j$ the set of polynomials of degree $j$ and $\delta_{ij}$ the Kronecker delta. One then obtains:
\begin{equation}\label{expf}
\dot y(ch) = f(y(ch)) \equiv \sum_{j\ge 0} P_j(c)\gamma_j(y), \qquad c\in[0,1],
\end{equation}
with the Fourier coefficients $\gamma_j(y)$ given by
\begin{equation}\label{gammaj}
\gamma_j(y)=\int_0^1 P_j(\tau)f(y(\tau h))\dd\tau, \quad j=0,1,\dots.
\end{equation}
We recall that:
\begin{equation}\label{Pjprop}
\|P_j\| := \max_{x\in[0,1]} |P_j(x)| = \sqrt{2j+1}, \qquad \int_0^1 P_j(x) q(x)\dd x = 0, \quad \forall q\in\Pi_{j-1}.
\end{equation}

Let us now study the properties of the coefficients $\gamma_j(y)$ defined at (\ref{gammaj}). To begin with, we report a result adapted from \cite[Lemma\,1]{BIT2012}.

\begin{lem}\label{hj} Let $g:[0,h]\rightarrow V$, with $V$ a vector space, admit a Taylor expansion at 0. Then
$$\int_0^1P_j(c) g(ch)\dd c =O(h^j).$$
\end{lem}
\proof Taking into account the second formula in (\ref{Pjprop}), we get
\begin{eqnarray*}
\int_0^1 P_j(c) g(c h)\dd c &=&\int_0^1 P_j(c) \sum_{\ell\ge0} \frac{g^{(\ell)}(0)}{\ell !} (c h)^\ell\dd c
~=~\sum_{\ell\ge j} \frac{g^{(\ell)}(0)}{\ell !} h^\ell \int_0^1 P_j(c) c^\ell\dd c\\
&\equiv& \sum_{\ell\ge j}g^{(\ell)}(0)  \phi_{j\ell} h^\ell ~=~O(h^j), 
\end{eqnarray*}
with ~$\phi_{j\ell} = \frac{1}{\ell !} \int_0^1 P_j(c) c^\ell\dd c$,~ so that (see the left formula in (\ref{Pjprop}))
$$|\phi_{j\ell}| \le \frac{\|P_j\|}{\ell!} = \frac{\sqrt{2j+1}}{\ell!},\qquad \ell=j,j+1,\dots.\,\QED$$\smallskip

\begin{cor}\label{Ohj} The Fourier coefficients defined in (\ref{gammaj}) satisfy: $\gamma_j(y) = O(h^j)$.\end{cor}

\medskip We now want to derive an estimate which generalizes the result of Corollary~\ref{Ohj} to the case where the stepsize $h$ is not small. For this purpose, hereafter we assume that the solution $y(t)$ of (\ref{ivp}) admit a  complex analytic extension in a neighbourhood of 0. Moreover, we shall denote by $\B(0,r)$ the closed ball of center 0 and radius $r$ in the complex plane, and $\C(0,r)$ the corresponding circumference. The following results then hold true. 

\begin{lem}\label{Davislem}
Let $P_j$ be the $j$th shifted and scaled Legendre polynomial and, for $\rho>1$, let us define the function  
\begin{equation}\label{Qj}
Q_j(\xi) =   \int_0^1 \frac{P_j(c)}{\xi-c}\dd c, \qquad \xi\in\C(0,\rho).
\end{equation}
Then, 
\begin{equation}\label{normQj}
\|Q_j\|_\rho := \max_{\xi\in\C(0,\rho)}|Q_j(\xi)| \le \frac{\sqrt{2j+1}}{(\rho-1)\rho^j}.
\end{equation}
\end{lem}
\proof One has, for $|\xi|=\rho>1$, and taking into account (\ref{Pjprop}):
\begin{eqnarray*}
Q_j(\xi) &=&  \int_0^1 \frac{P_j(c)}{\xi-c}\dd c   ~=~ \xi^{-1}  \int_0^1 \frac{P_j(c)}{1-\xi^{-1}c}\dd c ~=~
\xi^{-1}  \int_0^1 P_j(c)\sum_{\ell\ge0}\xi^{-\ell}c^\ell \,\dd c\\
&=& \xi^{-1} \sum_{\ell\ge j} \xi^{-\ell} \int_0^1 P_j(c) c^\ell \,\dd c ~=~ \xi^{-j-1} \sum_{\ell\ge0}\xi^{-\ell}  \int_0^1 P_j(c)c^{\ell+j} \,\dd c.
\end{eqnarray*}
Passing to norms, one has:
\begin{itemize}
\item $\left|\int_0^1 P_j(c)c^{\ell+j} \,\dd c\right|\le \|P_j\|=\sqrt{2j+1}$\,,
\item $\left| \xi^{-j-1} \sum_{\ell\ge0}\xi^{-\ell}\right| \le \rho^{-j-1} \sum_{\ell\ge0}\rho^{-\ell} = \left[ (\rho-1)\rho^j\right]^{-1}$\,,
\end{itemize}
from which (\ref{normQj}) follows.\,\QED
\bigskip

\begin{lem}\label{roh}
Let $g(z)$ be analytic in the closed ball $\B(0,r^*)$ of the complex plane, for a given $r^*>0$. Then,  for all $0<h<r^*$,  
\begin{equation}\label{ghz}
g_h(\xi) := g(\xi h)
\end{equation}
is analytic in $\B(0,\rho)$, with 
\begin{equation}\label{rho1}
\rho \,\equiv\, \rho(h) \,:=\,  \frac{r^*}h >1.
\end{equation}
\end{lem}
\bigskip

We are now in the position of stating the following result.\footnote{The used arguments are mainly adapted from \cite{Davis1975}.}

\smallskip 
\begin{theo}\label{Davisth} Assume that the function 
\begin{equation}\label{gz}
g(z):= f(y(z)) 
\end{equation}
satisfies the hypotheses of Lemma~\ref{roh}, and that $h>0$ in (\ref{expf})--(\ref{gammaj}) is such that (\ref{rho1}) holds true. Then, there exists $\kappa>0$, independent of $h$, such that\,
\footnote{Hereafter, for sake of clarity, we shall denote by $|\cdot|$ any convenient vector or matrix norm.}
\begin{equation}\label{gammajbound}
|\gamma_j(y)| \le \kappa\,\sqrt{2j+1}\,\rho^{-j}.
\end{equation}
\end{theo}
\proof By considering the function (\ref{ghz}) corresponding to (\ref{gz}), and with reference to the function $Q_j(\xi)$ defined in (\ref{Qj}), one has that the parameter $\rho$, as defined in (\ref{rho1}), is greater than 1 and, moreover, (see (\ref{gammaj}))
\begin{eqnarray*}
\gamma_j(y) &=& \int_0^1P_j(c)f(y(ch))\dd c ~\equiv~ \int_0^1 P_j(c) g_h(c)\dd c ~=~ \int_0^1 P_j(c)\left[\frac{1}{2\pi \ii}\int_{\C(0,\rho)} \frac{g_h(\xi)}{\xi-c}\dd \xi\right] \dd c \\[3mm]
&=&\frac{1}{2\pi \ii}\int_{\C(0,\rho)}g_h(\xi) \left[\int_0^1  \frac{P_j(c)}{\xi-c}\dd c\right] \dd \xi
~\equiv~\frac{1}{2\pi \ii}\int_{\C(0,\rho)}g_h(\xi)  Q_j(\xi)\, \dd \xi.
\end{eqnarray*}
Then, passing to norms (see (\ref{normQj})), 
$$|\gamma_j(y)|\le \rho\|g_h\|_\rho \|Q_j\|_\rho.$$
Moreover, observing that (see (\ref{gz}), (\ref{ghz}), and (\ref{rho1})):
$$\|g_h\|_\rho :=   \max_{\xi\in \C(0,\rho)} \left|g_h\left(\xi\right)\right|  \le \max_{\xi\in \B(0,\rho)} \left|g_h\left(\xi\right)\right| \equiv \max_{z\in \B(0,r^*)} |g(z)|=:\|g\|,$$
and using (\ref{normQj}), one has:
$$|\gamma_j(y)|\le  \frac{\rho}{(\rho-1)} \|g\| \sqrt{2j+1} \rho^{-j}.$$
Furthermore, by considering that for all $\rho\ge\rho^*>1$ 
$$
 \frac{\rho}{(\rho-1)} \|g\| 
$$
 is bounded, (\ref{gammajbound}) eventually follows.\,\QED\bigskip

\begin{rem}\label{boundongamma}
It is worth mentioning that, in the bound (\ref{gammajbound}), the dependence on $h$ only concerns the parameter $\rho>1$, via the expression (\ref{rho1}), from which one infers that $\rho\sim h^{-1}$, for all\, $0<h<r^*$. This, in turn,  is consistent with the result of Corollary~\ref{Ohj}, when $h\rightarrow 0$.
\end{rem}

Let us now consider a polynomial approximation to (\ref{expf}),
\begin{equation}\label{sig1}
\dot\sigma(ch) \,=\, \sum_{j=0}^{s-1} P_j(c) \gamma_j(\sigma),\qquad c\in[0,1],
\end{equation}
where $\gamma_j(\sigma)$ is defined according to (\ref{gammaj}) by formally replacing $y$ by $\sigma$, i.e., 
\begin{equation}\label{gammaj1}
\gamma_j(\sigma)=\int_0^1 P_j(\tau)f(\sigma(\tau h))\dd\tau, \quad j=0,1,\dots,s-1.
\end{equation}
Integrating term by term (\ref{sig1}), and imposing the initial condition in (\ref{ivp}), provide us with the polynomial approximation of degree $s$:
\begin{equation}\label{sig}
\sigma(ch) = y_0 + h\sum_{j=0}^{s-1} \int_0^c P_j(x)\dd x \, \gamma_j(\sigma), \qquad c\in[0,1].
\end{equation}
We now want to assess the extent to which $\sigma(ch)$ approximates $y(ch)$, for $c\in[0,1]$. When $h\rightarrow 0$, it is known that $y(h)-\sigma(h)=O(h^{2s+1})$ (see, e.g., \cite{LIMbook2016,BI2018,BIT2012}). Nevertheless, we here discuss the approximation of $\sigma$ to $y$, in the interval $[0,h]$, when $h$ is finite and only assuming that the result of Theorem~\ref{Davisth} is valid. 
The following result then holds true.\footnote{The proof uses arguments similar to those of \cite[Theorem~4]{BIT2012}.}

\begin{theo}\label{spectral}
Let $y$ be the solution of (\ref{ivp}), $\sigma$ be defined according to (\ref{sig}), and assume that
$f(\sigma(z))$ is analytic in $\B(0,r^*)$, for a given $r^*>0$.
Then, for all $0<h<r^*$, there exist $M,\hM>0$, independent of $h$, and $\rho,\hrho>1$, $\rho,\hrho\sim h^{-1}$, such that:
\begin{itemize}
\item $|\sigma(ch)-y(ch)|\le chM(2s+1)\,\rho^{-s} + h.o.t.$, \quad $c\in[0,1]$\,;
\item $|\sigma(h)-y(h)|\le h \hM(2s+1)\,\hrho^{-2s} + h.o.t.$\,.
\end{itemize}
\end{theo}

\proof 
Let $y(t,\xi,\eta)$ denote the solution of the problem
$$\dot y = f(y), \quad t\ge\xi, \qquad y(\xi)=\eta,$$ and $\Phi(t,\xi)$ be the solution of the associated variational problem,
$$\dot\Phi(t,\xi) = f'(y(t,\xi,\eta))\Phi(t,\xi), \quad t\ge\xi,\qquad \Phi(\xi,\xi)=I,$$ 
having set $f'$ the Jacobian of $f$. We also recall the following well-known perturbation results:
$$\frac{\partial}{\partial\eta} y(t,\xi,\eta) = \Phi(t,\xi), \qquad \frac{\partial}{\partial\xi} y(t,\xi,\eta) = -\Phi(t,\xi)f(\eta).$$
Consequently, from (\ref{gammaj1}) and (\ref{sig}), one has:

\begin{eqnarray}\nonumber
\sigma(ch)-y(ch)&=& y(ch,ch,\sigma(ch))-y(ch,0,\sigma(0)) 
~=~ \int_0^{ch} \frac{\dd}{\dd t} y(ch,t,\sigma(t))\,\dd t\\ \nonumber
&=&\int_0^{ch} \left[\left.\frac{\partial}{\partial\xi} y(ch,\xi,\sigma(t))\right|_{\xi=t}+
\left.\frac{\partial}{\partial\eta} y(ch,t,\eta)\right|_{\eta=\sigma(t)}\dot\sigma(t)\right]\dd t\\ \nonumber
&=&-h\int_0^c\Phi(ch,\tau h)\left[ f(\sigma(\tau h))-\dot\sigma(\tau h)\right]\dd \tau\\ \nonumber
&=& -h \int_0^c\Phi(ch,\tau h)\left[ \sum_{j\ge 0} P_j(\tau) \gamma_j(\sigma)-\sum_{j=0}^{s-1}P_j(\tau)\gamma_j(\sigma)\right]\dd \tau\\ \label{bound}
&=&-h\sum_{j\ge s}\left[\int_0^c P_j(\tau)\Phi(ch,\tau h)\dd\tau\right] \gamma_j(\sigma).
\end{eqnarray}
From the result of Theorem~\ref{Davisth} applied to $g(z):=f(\sigma(z))$, we know that there exist $\kappa$ independent of $h$ and $\rho>1$, $\rho\sim h^{-1}$, such that, for the Fourier coefficients defined in (\ref{gammaj1}),
\begin{equation}\label{rho}
|\gamma_j(\sigma)|\le\kappa\,\sqrt{2j+1}\,\rho^{-j}.
\end{equation}
Moreover, (see (\ref{Pjprop})) $\|P_j\| = \sqrt{2j+1}$ and, considering that, for all $h\in(0,r^*)$,
$$\max_{x_1,x_2\in[0,h]}|\Phi(x_1,x_2)| \le \max_{x_1,x_2\in[0,r^*]}|\Phi(x_1,x_2)|=:\nu,$$
which is independent of $h$, the first statement then follows from (\ref{bound}), by setting $M=\nu\kappa$, since
$$|\sigma(ch)-y(ch)|\le c h M \sum_{j\ge s} \sqrt{2j+1}\,\rho^{-j},\qquad c\in[0,1].$$
To prove the second statement (i.e., when $c=1$), we observe that the result of Theorem~\ref{Davisth} holds true also for the Fourier coefficients
$$
\int_0^1 P_j(\tau)\Phi(ch,\tau h)\dd\tau$$ by setting $g(z):=\Phi(ch,z)$. Consequently, there exist $\kappa_1>0$, independent of $h$, and $\rho_1>1,$ $\rho_1\sim h^{-1}$, such that
\begin{equation}\label{rho_1}
\left|\int_0^1 P_j(\tau)\Phi(ch,\tau h)\dd\tau\right|\le \kappa_1\,\sqrt{2j+1}\,\rho_1^{-j}.
\end{equation}
The second statement then follows again from (\ref{bound}) by setting \,$\hM=\kappa_1\kappa$\,
and \,$\hrho=\max\{\rho_1,\rho\}$, so that:
$$|\sigma(h)-y(h)| \le h \hM \sum_{j\ge s} (2j+1)\hrho^{\,-2j}.\,\QED$$

Let us now introduce the use of a finite precision arithmetic, with machine precision $u$, for approximating (\ref{expf}). Then, the best we can do is to consider the polynomial approximation (\ref{sig1})--(\ref{gammaj1})\,\footnote{Hereafter, $\doteq$ means ``equal within the round-off error level of the used finite precision arithmetic''.}
\begin{equation}\label{sig1u}
\dot y(ch) \,\doteq\, \dot\sigma(ch) \,=\, \sum_{j=0}^{s-1} P_j(c) \gamma_j(\sigma),\qquad c\in[0,1],
\end{equation}
such that
\begin{equation}\label{gammas}
|\gamma_s(\sigma)| < tol \cdot \max_{j<s} |\gamma_j(\sigma)|, \qquad tol\sim u.
\end{equation}
Integrating (\ref{sig1u}), and imposing that $\sigma(0)=y_0$, then brings back to (\ref{sig}). 
We observe that because of (\ref{rho}), (\ref{gammas})  may be approximately recast as
\begin{equation}\label{gammas1}
\sqrt{2s+1}\,\rho^{-s} < tol \sim u,
\end{equation}
where $\rho\sim h^{-1}$. Consequently, choosing $s$ such that (\ref{gammas}) (or (\ref{gammas1})) is satisfied, we obtain that:
\begin{itemize}
\item the polynomial $\sigma(ch)$ defined by (\ref{sig1u}) and (\ref{sig}) provides a uniformly accurate approximation to $y(ch)$, in the whole interval $[0,h]$, within the possibility of the used finite precision arithmetic;

\item $\sigma(h)$ is a {\em spectrally accurate} approximation to $y(h)$. In particular, in light of the second point of the result of Theorem~\ref{spectral}, one has that the criterion (\ref{gammas}) can be conveniently relaxed. In fact,  making the ansatz (see (\ref{rho}) and (\ref{rho_1})) $\rho=\rho_1$ and $\kappa=\kappa_1$, one has that
\begin{equation}\label{superco}
|\sigma(h)-y(h)| \lesssim h\kappa^2 (2s+1)\rho^{-2s} \approx h |\gamma_s(\sigma)|^2.
\end{equation}
Imposing the approximate upper bound to be smaller than the machine epsilon $u$, one then obtains:
\begin{equation}\label{gammas2}
|\gamma_s(\sigma)|~\lesssim~ \sqrt{\frac{u}h} ~\propto~ u^{1/2},
\end{equation}
which is generally much less restrictive than (\ref{gammas}).\footnote{This latter criterion was that used in \cite{BMR2018} and \cite{BIMR2018}.} Alternatively, by considering that the use of relatively large time-steps $h$ is sought, one can use $tol\sim u^{1/2}$ in (\ref{gammas}), that is, 
\begin{equation}\label{newtol}
|\gamma_s(\sigma)| < tol \cdot \max_{j<s} |\gamma_j(\sigma)|, \qquad tol\sim u^{1/2},
\end{equation}
In other words, (\ref{superco}) means that the method maintains the property of {\em super-convergence}, which is known to hold when $h\rightarrow0$, also in the case where the time-step $h$ is relatively large.
\end{itemize}
\begin{rem}
In particular, we observe that (\ref{gammas}) (or (\ref{gammas1}) or (\ref{gammas2})) can be fulfilled by varying
the value of $s$, and/or that of the stepsize $h$, by considering that, by virtue of (\ref{rho1}),
$$
\rho(h_{new}) \approx \rho(h_{old})\frac{h_{old}}{h_{new}},
$$
$h_{old}$ and $h_{new}$ being the old and new stepsizes, respectively.

It is worth mentioning that the result of Theorem~\ref{spectral} can be also used to define a stepsize variation, within a generic error tolerance $tol$, thus defining a strategy for the simultaneous order/stepsize variation.

\end{rem}

We conclude this section mentioning that, to gain efficiency, the criterion (\ref{gammas}) for the choice of $s$ in (\ref{sig1u}) can be more conveniently changed to
\begin{equation}\label{tol}
|\gamma_{s-1}(\sigma)| \le tol\,\cdot \max_{j\le s-1} |\gamma_j(\sigma)|, \qquad tol \sim \rho\cdot u.
\end{equation}
Similarly, the less restrictive criterion (\ref{gammas2}) becomes
$$
|\gamma_{s-1}(\sigma)|\lesssim \rho\sqrt{\frac{u}h} \propto u^{1/2},
$$
or, alternatively, one uses  $tol\sim u^{1/2}$ in (\ref{tol}). 
As is clear, computing the norms of the coefficients $\gamma_j(\sigma)$ permits to derive estimates for the parameters 
$\kappa$ and $\rho$ in (\ref{rho}), as we shall see later in the numerical tests.

\section{SHBVMs}\label{shbvms}

The approximation procedure studied in the previous section does not yet provide a numerical method, in that the integrals defining $\gamma_j(\sigma)$, $j=0,\dots,s-1$,  in (\ref{gammaj1})--(\ref{sig}) need to be computed. For this purpose,  one can approximate them to within machine precision  through a Gauss-Legendre quadrature formula of order $2k$ (i.e., the interpolatory quadrature rule defined at the zeros of $P_k$) with $k$ large enough. In particular, following the criterion used in \cite{BMR2018,BIMR2018}, for the double precision IEEE\,\footnote{In such a case, the machine precision is $u\approx 10^{-16}$.} we choose
\begin{equation}\label{kappas} k = \max\{ 20, s+2\}.\end{equation}
After that, we define the approximation to $y(h)$ as 
\begin{equation}\label{y1}y_1 :=\sigma(h) \equiv y_0+h\gamma_0(\sigma).\end{equation} 
In so doing, one eventually obtains a HBVM$(k,s)$ method, which we sketch below. Hereafter, we shall refer to such a method as to {\em spectral HBVM (in short, SHBVM)}, since its parameters $s$ and $k$, respectively defined in (\ref{gammas}) (or (\ref{gammas1}) or (\ref{gammas2})) and (\ref{kappas}), are aimed at obtaining a numerical solution which is accurate within the round-off error level of the used finite precision arithmetic.

For sake of completeness, let us now briefly sketch what a HBVM$(k,s)$ method is. In general, to approximate the Fourier coefficient $\gamma_j(\sigma)$, and assuming for sake of simplicity that full machine accuracy is gained, we use the quadrature
\begin{equation}\label{hgj}
\gamma_j(\sigma) \doteq \sum_{\ell=1}^k b_\ell P_j(c_\ell)f(\sigma(c_\ell h)) ~=:~\hg_j,\qquad j=0,\dots,s-1,
\end{equation}
where the polynomial $\sigma$ is that defined is (\ref{sig}) by formally replacing $\gamma_j(\sigma)$ with $\hg_j$, and $(c_i,b_i)$ are the abscissae and weights of the Gauss-Legendre quadrature of order $2k$ on the interval $[0,1]$.\footnote{I.e., $0<c_1<\dots<c_k<1$ are the zeros of $P_k$.}
Setting $Y_\ell=\sigma(c_\ell h)$, from (\ref{hgj}) one then obtains the {\em stage equations}
\begin{equation}\label{RK1}
Y_i ~=~ y_0 + h\sum_{j=0}^{s-1} \int_0^{c_i} P_j(x)\dd x \hg_j ~\equiv~ y_0 + h\sum_{j=1}^k \underbrace{b_j \left[\sum_{\ell=0}^{s-1} \int_0^{c_i}P_\ell(x)\dd x P_\ell(c_j)\right]}_{=:\,a_{ij}} f(Y_j), \quad i=1,\dots,k,
\end{equation}
with the new approximation given by (see (\ref{y1}))
\begin{equation}\label{RK2}
y_1 ~=~ y_0+ h\hg_0 ~\equiv~ y_0+h\sum_{i=1}^k b_i f(Y_i).
\end{equation}
Consequently, with reference to (\ref{RK1}), setting
\begin{equation}\label{Abc}
A = \left( a_{ij}\right)\in\RR^{k\times k}, \qquad \bfb = (b_i), ~ \bfc=(c_i)\in\RR^k,
\end{equation}
one easily realizes that (\ref{RK1}) and (\ref{RK2}) define the $k$-stage Runge-Kutta method with Butcher tableau:
$$ \begin{array}{c|c} \bfc & A\\ \hline \\[-3mm] &\bfb^\top\end{array}~. $$
From (\ref{RK1}) one verifies that the Butcher matrix in (\ref{Abc}) can be written as
\begin{equation}\label{A}A = \I_s\P_s^\top\Omega,  \end{equation}
with
\begin{equation}\label{PI}
\P_s = \pmatrix{ccc} P_0(c_1)&\dots&P_{s-1}(c_1)\\ \vdots & &\vdots\\ P_0(c_k)&\dots&P_{s-1}(c_k)\endpmatrix,~
\I_s = \pmatrix{ccc} \int_0^{c_1}P_0(x)\dd x&\dots&\int_0^{c_1}P_{s-1}(x)\dd x\\ \vdots & &\vdots\\ \int_0^{c_k}P_0(x)\dd x&\dots&\int_0^{c_k}P_{s-1}(x)\dd x\endpmatrix\in\RR^{k\times s},
\end{equation}
and
\begin{equation}\label{Om}
\Omega = \pmatrix{ccc} b_1\\ &\ddots\\ &&b_k\endpmatrix\in\RR^{k\times k}.
\end{equation}
In fact, setting $\bfe_i\in\RR^k$ the $i$th unit vector, and taking into account (\ref{A})--(\ref{Om}), one has
$$
\bfe_i^\top A\bfe_j =\bfe_i^\top \I_s \P_s^\top \Omega\bfe_j = \bfe_i^\top \I_s \left( \bfe_j^\top \P_s \right)^\top b_j =b_j \left[\sum_{\ell=0}^{s-1} \int_0^{c_i}P_\ell(x)\dd x P_\ell(c_j)\right] \equiv a_{ij},
$$
as defined in (\ref{RK1}). 
For later use, we also recall the following known property relating the matrices (\ref{PI})--(\ref{Om}) (see, e.g., \cite[Lemma~3.6]{LIMbook2016}):
\begin{equation}\label{Xs}
\P_s^\top\Omega\I_s = \pmatrix{cccc} \xi_0 & -\xi_1\\ \xi_1 &0 &\ddots\\ &\ddots & \ddots &-\xi_{s-1}\\ & &\xi_{s-1} &0\endpmatrix =: X_s\in\RR^{s\times s}, \qquad \xi_i = \left(2\sqrt{ |4i^2-1|}\right)^{-1}.
\end{equation}

At this point, we observe that the stage equations (\ref{RK1}) can be cast in vector form, by taking into account (\ref{Abc})--(\ref{Om}), as
\begin{equation}\label{Y}
Y \equiv \pmatrix{c} Y_1\\ \vdots \\ Y_k\endpmatrix = \bfe\otimes y_0 + h\I_sP_s^\top\Omega \otimes I_m \cdot f(Y), \qquad \bfe=\pmatrix{c}1\\ \vdots \\1\endpmatrix\in\RR^k,
\end{equation}
with an obvious meaning of $f(Y)$. On the other hand, the block vector of the coefficients in (\ref{hgj}) turns out to be given by
\begin{equation}\label{hbg}
\hbfg \equiv \pmatrix{c} \hg_0 \\ \vdots \\ \hg_{s-1}\endpmatrix = \P_s^\top\Omega\otimes I_m\cdot f(Y).
\end{equation}
Consequently, from (\ref{Y}) one obtains
$$Y = \bfe\otimes y_0 + h\I_s\otimes I_m\cdot \hbfg,$$
and then, from (\ref{hbg}), one eventually derives the equivalent discrete problem
\begin{equation}\label{F}
F(\hbfg) \,:=\, \hbfg - \P_s^\top\Omega\otimes I_m\cdot f\left(\bfe\otimes y_0 + h\I_s\otimes I_m\cdot \hbfg\right) \,=\, {\bf0},
\end{equation}
which has (block) dimension $s$, independently of $k$ (compare with (\ref{Y})). Once it has been solved, the new approximation is obtained (see (\ref{RK2})) as $y_1=y_0+h\hg_0$. 

It is worth observing that the new discrete problem (\ref{F}), having block dimension $s$ independently of $k$, allows us to use arbitrarily high-order quadratures (see (\ref{kappas})), without affecting that much the computational cost.

In order to solve (\ref{F}), one could in principle use a fixed-point iteration,\footnote{Hereafter, the initial approximation $\hbfg^0=\bf0$ is conveniently used.}
$$ \hbfg^{\ell+1} :=  \P_s^\top\Omega\otimes I_m\cdot f\left(\bfe\otimes y_0 + h\I_s\otimes I_m\cdot \hbfg^\ell\right), \qquad \ell=0,1,\dots,$$
which, though straightforward, usually implies restrictions on the choice of the  stepsize $h$. For this reason, this approach is generally not useful when using the methods as spectral methods, where the use of relatively large stepsizes is sought. 
On the other hand, the use of the simplified Newton iteration for solving (\ref{F})  reads, by virtue of (\ref{Xs}),
\begin{equation}\label{Newtsim}
\mbox{solve:\,} \left[ I_s\otimes I_m-hX_s\otimes f'(y_0)\right] \bfdelta^\ell = -F(\hbfg^\ell), \qquad   \hbfg^{\ell+1} := \hbfg^\ell+\bfdelta^\ell, \qquad \ell=0,1,\dots.
\end{equation}
However, the coefficient matrix in (\ref{Newtsim}) has a dimension $s$ times larger than that of the continuous problem (i.e., $m$) and, therefore, this can be an issue when large value of $s$ are to be used, as in the case of SHBVMs. Fortunately, this problem can be overcome by replacing the previous iteration (\ref{Newtsim}) with a corresponding {\em blended iteration} \cite{BIT2011,LIMbook2016,BI2018} (see also \cite{BFCI2014}). In more details, once one has formally computed the $m\times m$ matrix
\begin{equation}\label{blend1}
\Sigma = \left(I_m-h\rho_s f'(y_0)\right)^{-1}, \qquad \rho_s = \min_{\lambda\in\sigma(X_s)}|\lambda|,
\end{equation}
one iterates:
\begin{equation}\label{blend2}
\bfeta^\ell:=F(\hbfg^\ell),\quad \bfeta_1^\ell := \rho_sX_s^{-1}\otimes I_m\bfeta^\ell, \quad 
\hbfg^{\ell+1} := \hbfg^\ell +I_s\otimes \Sigma \left[ \bfeta_1^\ell + I_s\otimes \Sigma\left( \bfeta^\ell-\bfeta_1^\ell\right)\right], \qquad \ell=0,1,\dots.
\end{equation}
Consequently, one only needs to compute, at each time-step, the matrix $\Sigma\in\RR^{m\times m}$ defined in (\ref{blend1}),\footnote{I.e., factor $\Sigma^{-1}$.} having the same size as that of the continuous problem. Moreover, it is worth mentioning that for  semi-linear problems with a leading linear part, the Jacobian of $f$ can be approximated with the (constant) linear part, so that $\Sigma$ is computed once for all \cite{BFCI2015,BBFCI2018,BGS2018,BZL2018,BMR2018,BIMR2018}.

\begin{rem}\label{sticazzi}
It must be stressed that it is the availability of the very efficient blended iteration (\ref{blend1})--(\ref{blend2}) which makes possible the practical use of HBVMs as spectral methods in time, since relatively large values of $s$ can be easily and effectively handled. 

A thorough analysis of the blended iteration can be found in \cite{BrMa2002}. Contexts where it has  been successfully implemented include stiff ODE-IVPs \cite{BrMa2004}, linearly implicit DAEs up to index 3 \cite{BrMa2006} (see also the  code {\tt BiMD} in TestSet for IVP Solvers \cite{TSet}), and canonical Hamiltonian systems (see the Matlab code {\tt HBVM}, available at \cite{LIMwebpage}), while its implementation in the solution of RKN methods may be found in \cite{WMF2017}.
\end{rem}

\section{Numerical tests}\label{numtest}
The aim of this section is twofold: firstly, to assess the theoretical analysis of SHBVMs made in Section~\ref{spectime}; 
secondly, to compare such methods w.r.t. some well-known ones.   All numerical tests, which concern different kinds of ODE problems, have been computed on a laptop with a 2.8GHz Intel-i7 quad-core processor and 16GB of memory, running Matlab 2017b.  For the SHBVM method, the criterions (\ref{newtol}) and (\ref{kappas}) have been respectively used to determine its parameters $s$ and $\kappa$.
 
\subsubsection*{The Kepler problem}
We start considering the well-known Kepler problem (see, e.g., \cite[Chapter~2.5]{LIMbook2016}), which is Hamiltonian, with Hamiltonian function
\begin{equation}\label{keplH}
H(q,p) = \frac{1}2 \|p\|^2 - \|q\|_2^{-1}, \qquad q,p\in\RR^2.
\end{equation}
Consequently, we obtain the equations
\begin{equation}\label{kepl}
\dot q = p, \qquad \dot p = -\|q\|_2^{-3}q,
\end{equation}
which, when coupled with the initial conditions
\begin{equation}\label{kepl0}
q(0) = \pmatrix{cc}1-\eps,&0\endpmatrix^\top, \qquad p(0) =  \pmatrix{cc} 0, & \sqrt{\frac{1+\eps}{1-\eps}}\endpmatrix^\top, \qquad \eps\in[0,1),
\end{equation}
provide a periodic orbit of period $T=2\pi$ that, in the $q$-plane, is given by an ellipse of eccentricity $\eps$. In particular, we choose the value $\eps=0.5$. The solution of this problem has two additional (functionally independent) invariants besides the Hamiltonian (\ref{keplH}), i.e., the {\em angular momentum} and one of the nonzero components of the Lenz vector \cite[page\,64]{LIMbook2016} (in particular, we select the second one):
\begin{equation}\label{keplML}
M(q,p) = q_1p_2-p_1q_2, \qquad L(q,p) = -p_1M(q,p)-q_2\|q\|_2^{-1}.
\end{equation}

\begin{figure}[t]
\centerline{\includegraphics[width=12cm,height=9cm]{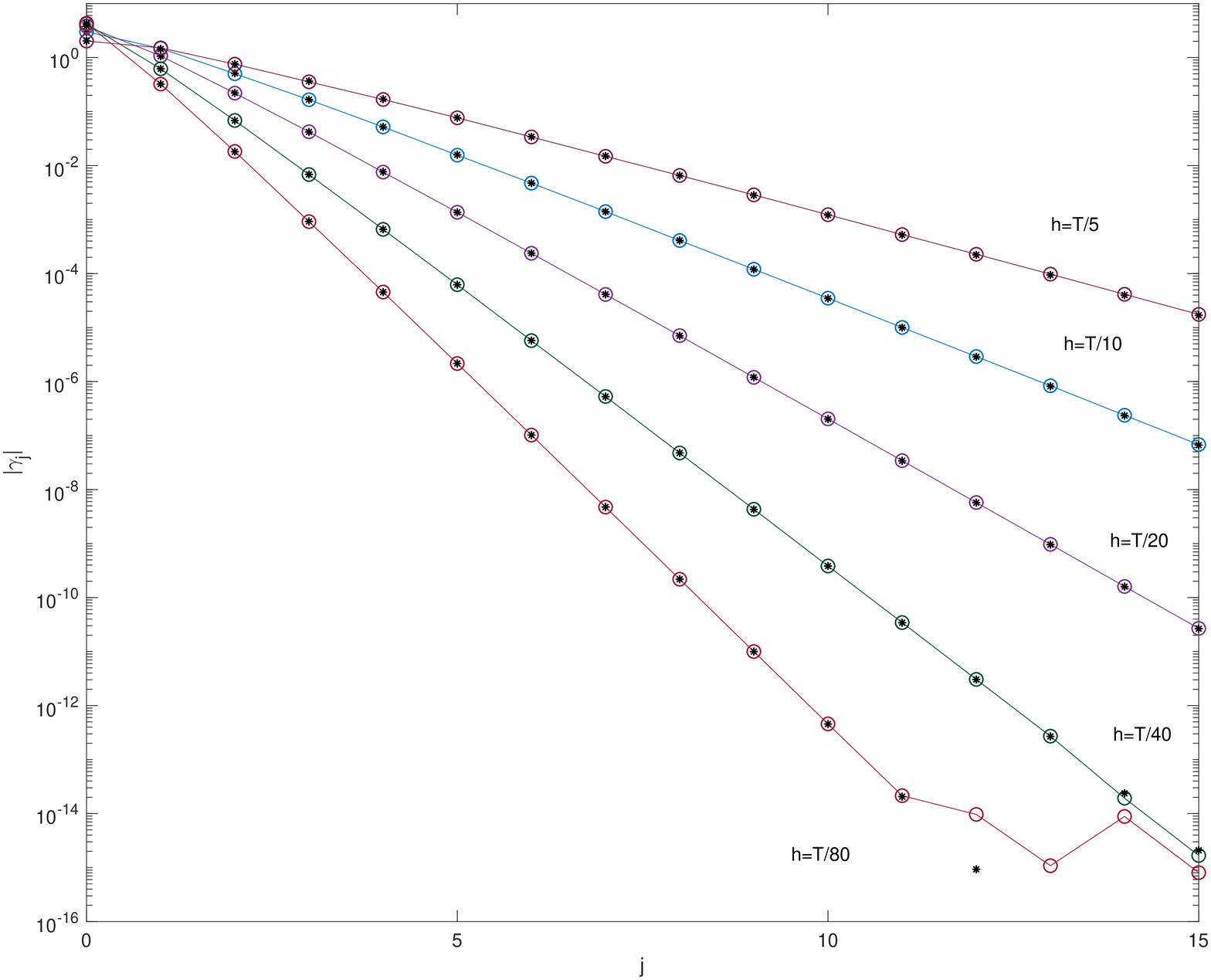}}
\caption{behavior of $|\hat\gamma_j|$ for decreasing values of the time-step $h$ for the Kepler problem (\ref{kepl})--(\ref{kepl0}) solved by the HBVM(20,16) method with decreasing time-steps. The line with circles are the computed norms, whereas the asterisks are the estimated ones. Observe that, for the smallest time-steps, the computed norms stagnate near the round-off error level.}
\label{rok}
\end{figure}

At first, we want to assess the result of Theorem~\ref{Davisth}. For this purpose, we apply the HBVM(20,16) method for one step starting from the initial condition (\ref{kepl0}), and using time-steps $h_i = 2\pi/(5\cdot 2^{i-1})$, $i=1,\dots,5$. In Figure~\ref{rok} is the plot (see (\ref{hgj})) of $|\hat\gamma_j|$, for $j=0,1,\dots,15$, (solid line with circles), which, according to (\ref{rho}), should behave as $\kappa\sqrt{2j+1}\rho^{-j}$, due to the result of Theorem~\ref{Davisth}. A least square approximation technique has been employed to estimate the two parameters $\kappa$ and $\rho$  appearing in the bound (\ref{rho}). These theoretical bounds are highlighted by asterisks in Figure \ref{rok}. On the other hand, the line with circles in the figure are the computed norms of the vectors $\hat\gamma_j$: evidently, they well fit the theoretical estimations, except those which are close to the round-off error level. Moreover, according to the arguments in the proof of Theorem~\ref{Davisth}, one also expects that the estimate of $\kappa$ increases but is bounded, as $h\rightarrow 0$, whereas $\rho$ should be proportional to $h^{-1}$. This fact is confirmed by the results listed in Table~\ref{kro}. 

Next, we compare the following methods for solving (\ref{kepl})--(\ref{kepl0}):

\begin{table}[t]
\caption{estimated values for the parameters $\kappa$ and $\rho$ for Kepler problem (\ref{kepl})--(\ref{kepl0}), when using decreasing time-steps.}
\label{kro}
\smallskip
\centerline{\begin{tabular}{|lcr|}
\hline
$h$      & $\kappa$ & $\rho$\\
\hline
$2\pi/5$   & 2.0 & 2.4 \\
$2\pi/10$ & 3.0 & 3.6 \\
$2\pi/20$ & 3.8 & 6.2 \\
$2\pi/40$ & 4.2 & 11.7\\
$2\pi/80$ & 4.3 & 23.1\\
\hline\end{tabular}}
\end{table}

\begin{itemize}
\item the $s$-stage Gauss method (i.e., HBVM$(s,s)$), $s=1,2$, which is symplectic and of order $2s$. Consequently, it is expected to conserve the angular momentum $M(q,p)$ in (\ref{keplML}), which is a quadratic invariant;
\item the HBVM(6,$s$) method, $s=1,2$, which, for the considered stepsizes is energy-conserving and of order $2s$;
\item the SHBVM method described above, where $s$ and $k$ are determined according to (\ref{newtol}) and (\ref{kappas}), respectively, with $tol\approx 10^{-8}$. This tolerance, in turn, should provide us with full accuracy, according to the result of Theorem~\ref{spectral}, because of the super-convergence of the method, which is valid for any used step-size.\footnote{As matter of fact, considering more stringent tolerances does not improve the accuracy of the computed numerical solution.}
\end{itemize}
It is worth mentioning that the execution times  that we shall list for the Gauss, HBVM, and SHBVM methods are perfectly comparable, since the same Matlab code has been used for all such methods. This code, in turn, is a slight variant of the {\tt hbvm } function available at the url \cite{LIMwebpage}.

In Tables~\ref{kepl1}--\ref{kepl3}, we list the obtained results when using a time-step $h=2\pi/n$ over 100 periods. In more details, we list the maximum errors,  measured at each period, in the invariants (\ref{keplH}) and (\ref{keplML}), $e_H,e_M,e_L$, respectively, the solution error, $e_y$, and the execution times (in sec). As is expected, the symplectic methods conserve the angular momentum (since it is a quadratic invariant), whereas the energy-conserving HBVMs conserve the Hamiltonian function.\footnote{In this case, the Gauss methods exhibit a super-convergence in the conservation of the Hamiltonian (3 times the usual order) and HBVMs do the same with the angular momentum. This is due to the fact that the error is measured only at the end of each period.}  On the other hand, the SHBVM conserves all the invariants, and has a uniformly small solution error, by using very large stepsizes. Further, its execution time is the lowest one (less than 0.5 sec, when using $h=2\pi/5$), thus confirming the effectiveness of the method.

\begin{table}[ph]
\small
\caption{numerical result for the $s$-stage Gauss method, $s=1,2$, used for solving the Kepler problem (\ref{kepl})--(\ref{kepl0}), $\eps=0.5$, with stepsize $h=2\pi/n$.}
\label{kepl1}
\centerline{\begin{tabular}{|r|r|c|c|c|c|c|c|c|}
\hline
\multicolumn{9}{|c|}{Gauss-1}\\
\hline
$n$ & time & $e_H$ & rate & $e_M$ & $e_L$ & rate &$e_y$ & rate\\
\hline
 100 &   2.52 & 6.56e-03 & --- & 5.88e-15 & 4.97e-01 & --- & 3.04e\,00 & --- \\ 
  200 &   4.74 & 1.63e-03 & 2.0 & 1.04e-14 & 3.54e-01 & 0.5 & 2.39e\,00 & 0.3 \\ 
  400 &   9.78 & 3.82e-04 & 2.1 & 2.09e-14 & 9.76e-02 & 1.9 & 1.61e\,00 & 0.6 \\ 
  800 &  17.26 & 3.05e-05 & 3.6 & 7.66e-15 & 2.45e-02 & 2.0 & 7.49e-01 & 1.1 \\ 
 1600 &  33.74 & 6.07e-07 & 5.6 & 1.93e-14 & 6.11e-03 & 2.0 & 2.08e-01 & 1.8 \\ 
 3200 &  65.46 & 9.65e-09 & 6.0 & 3.04e-14 & 1.53e-03 & 2.0 & 5.25e-02 & 2.0 \\ 
\hline
\multicolumn{9}{|c|}{Gauss-2}\\
\hline
$n$ & time & $e_H$ & rate & $e_M$ & $e_L$ & rate &$e_y$ & rate\\
\hline
  50 &   2.33 & 2.05e-06 & --- & 3.44e-15 & 3.81e-02 & --- & 3.17e-01 & --- \\ 
 100 &   4.09 & 5.37e-10 & 11.9 & 5.77e-15 & 2.43e-03 & 4.0 & 2.09e-02 & 3.9 \\ 
 200 &   7.52 & 1.44e-13 & 11.9 & 7.55e-15 & 1.53e-04 & 4.0 & 1.32e-03 & 4.0 \\ 
 400 &  14.48 & 9.55e-15 & 3.9 & 9.99e-14 & 9.55e-06 & 4.0 & 8.29e-05 & 4.0 \\ 
 800 &  26.75 & 1.53e-14 & *** & 1.49e-14 & 5.97e-07 & 4.0 & 5.18e-06 & 4.0 \\ 
 1600 &  52.46 & 3.81e-14 & *** & 1.95e-14 & 3.73e-08 & 4.0 & 3.24e-07 & 4.0 \\ 
 3200 &  101.74 & 3.49e-14 & *** & 4.71e-14 & 2.33e-09 & 4.0 & 2.04e-08 & 4.0 \\ 
\hline
\end{tabular}}

\medskip
\caption{numerical result for the HBVM(6,$s$) method, $s=1,2$, used for solving the  Kepler problem (\ref{kepl})--(\ref{kepl0}), $\eps=0.5$, with stepsize $h=2\pi/n$.}
\label{kepl2}
\centerline{\begin{tabular}{|r|r|c|c|c|c|c|c|c|}
\hline
\multicolumn{9}{|c|}{HBVM(6,1)}\\
\hline
$n$ & time & $e_H$ & $e_M$  & rate & $e_L$ & rate &$e_y$ & rate\\
\hline
  100 &   3.55 & 4.44e-16 & 9.09e-04 & --- & 4.99e-01 & --- & 2.94e\,00 & --- \\ 
  200 &   7.10 & 4.44e-16 & 2.12e-05 & 6.0 & 3.52e-01 & 1.9 & 9.68e-01 & 1.9 \\ 
  400 &  12.47 & 6.66e-16 & 3.39e-07 & 6.0 & 9.70e-02 & 2.0 & 2.58e-01 & 2.0 \\ 
  800 &  22.86 & 4.44e-16 & 5.29e-09 & 6.0 & 2.44e-02 & 2.0 & 6.46e-02 & 2.0 \\ 
 1600 &  45.46 & 4.44e-16 & 8.26e-11 & 6.0 & 6.10e-03 & 2.0 & 1.62e-02 & 2.0 \\ 
 3200 &  86.34 & 6.66e-16 & 1.30e-12 & 6.0 & 1.53e-03 & 2.0 & 4.04e-03 & 2.0 \\ 
\hline
\multicolumn{9}{|c|}{HBVM(6,2)}\\
\hline
$n$ & time & $e_H$ & $e_M$  & rate & $e_L$ & rate &$e_y$ & rate\\
\hline
  50 &   2.92 & 4.44e-16 & 1.09e-07 & --- & 3.82e-02 & --- & 4.64e-02 & --- \\ 
 100 &   4.50 & 4.44e-16 & 2.72e-11 & 12.0 & 2.43e-03 & 4.0 & 2.94e-03 & 4.0 \\ 
 200 &   8.10 &  4.44e-16  & 5.88e-15 & 12.1 & 1.53e-04 & 4.0 & 1.84e-04 & 4.0 \\ 
 400 &  15.48 & 4.44e-16 & 3.89e-15 & *** & 9.55e-06 & 4.0 & 1.15e-05 & 4.0 \\ 
 800 &  28.42 & 4.44e-16 & 1.40e-14 & *** & 5.97e-07 & 4.0 & 7.20e-07 & 4.0 \\ 
 1600 &  52.29 & 6.66e-16 & 1.73e-14 & *** & 3.73e-08 & 4.0 & 4.50e-08 & 4.0 \\ 
 3200 & 107.41 & 6.66e-16 & 1.40e-14 & *** & 2.33e-09 & 4.0 & 2.81e-09 & 4.0 \\ 
\hline
\end{tabular}}

\medskip
\caption{numerical result for the SHBVM method used for solving the Kepler problem (\ref{kepl})--(\ref{kepl0}), $\eps=0.5$, with stepsize $h=2\pi/n$.}
\label{kepl3}
\centerline{\begin{tabular}{|r|r|r|c|c|c|c|c|}
\hline
$n$ & $k$ & $s$ & time & $e_H$ & $e_M$ & $e_L$ & $e_y$ \\
\hline
     5 & 24 & 22  &  0.47   & 4.44e-16 & 2.01e-14 & 1.66e-14 & 8.00e-13  \\
   10 & 20 & 16 &   0.71   & 4.44e-16 & 6.22e-15 & 2.34e-14 & 6.13e-13 \\ 
   20 & 20 & 11 &   1.22   & 4.44e-16 & 6.66e-16 & 3.89e-15 & 3.87e-13 \\ 
   40 & 20 &   9 &   2.16   & 2.22e-16 & 1.89e-15 & 3.28e-15 & 5.75e-13 \\ 
\hline
\end{tabular}}
\end{table}

\subsubsection*{A Lotka-Volterra problem}
We consider the following Poisson problem \cite{CH2011},
\begin{equation}\label{lotka1}
\dot y = B(y) \nabla H(y), \qquad  B(y)^\top = -B(y),
\end{equation}
with $ y\in\RR^3$, 
\begin{equation}\label{lotka2}
B(y) = \pmatrix{ccc} 
0 & c\,y_1y_2 &bc\,y_1y_3\\
-c\,y_1y_2 &0& -y_2y_3\\
-bc\,y_1y_3 & y_2y_3 &0
\endpmatrix, \qquad H(y) = ab\, y_1 + y_2 -a\,y_3 +\nu\ln y_2 -\mu\ln y_3,
\end{equation}
and $abc=-1$. In this case, there is a further invariant besides the Hamiltonian $H$, i.e., the {\em Casimir} 
\begin{equation}\label{casimir}
C(y) = ab\ln y_1 -b\ln y_2 +\ln y_3.
\end{equation}
The solution turns out to be periodic, with period $T\approx2.878130103817$, when choosing
\begin{equation}\label{lotka3}
a = -2, \quad b = -1, \quad c = -0.5, \quad \nu = 1, \quad\mu = 2, \quad
y(0) = (1,\, 1.9,\, 0.5)^\top.
\end{equation}
For this problem, the HBVM$(k,s)$ method is no more energy-conserving, as well as the $s$-stage Gauss method. As matter of fact, both exhibit a drift in the invariants and a quadratic error growth in the numerical solution.  The obtained results for the SHBVM method, with $tol\approx 10^{-8}$ in (\ref{newtol}) for choosing $s$, $\kappa$ given by (\ref{kappas}), and using a stepsize $h=T/n$, are listed in Table~\ref{lotka}, where it is reported the maximum Hamiltonian error, $e_H$, the Casimir error, $e_C$, and the solution error $e_y$, measured at each period,  over 100 periods. In such a case, all the invariants turn out to be numerically conserved, and the solution error is uniformly very small. Moreover, the SHBVM using the largest time-step (i.e., $h=T/5\approx 0.57$) turns out to be the most efficient one. For comparison, in the table we also list the results obtained by using the Matlab solver {\tt ode45} used with the default parameters, requiring 5600 integration steps and stepsizes approximately in the range $[2.2\cdot 10^{-2},1.1\cdot 10^{-1}]$, and the same solver used with parameters {\tt AbsTol=1e-15}, {\tt RelTol=1e-10}, now requiring 121760 integration steps, with stepsizes approximately in the range $[10^{-3}, 4.2\cdot 10^{-3}]$.

\begin{table}[t]
\caption{numerical result for the SHBVM method used for solving the Lotka-Volterra problem (\ref{lotka1})--(\ref{lotka3}) with stepsize $h=T/n$. We also list the results obtained by using {\tt ode45}, both with the default parameters, and with parameters {\tt AbsTol=1e-15}, {\tt RelTol=1e-10} (which we denote by {\tt ode45*}).}
\label{lotka}
\centerline{\begin{tabular}{|r|r|r|c|c|c|c|}
\hline
$n$ & $k$ & $s$ & time & $e_H$ & $e_C$ &$e_y$ \\
\hline
   5 & 20 & 16 &   0.88 & 8.26e-14 &  4.89e-14 &  4.24e-11  \\ 
  10 & 20 & 11 &   1.37 & 1.33e-14 & 1.33e-14 &  5.01e-11  \\ 
  15 & 20 & 9 &   1.84 & 3.11e-14 &  1.62e-14 &  4.92e-11  \\ 
\hline
\multicolumn{3}{|c|}{\tt ode45} & 0.23 & 7.41e-01& 7.27e-01 & 3.62e\,00\\
\multicolumn{3}{|c|}{\tt ode45*} & 4.12 &  1.14e-08&   8.71e-09 & 8.44e-07\\
\hline
\end{tabular}}

\medskip
\caption{numerical result for the SHBVM method used for solving the stiff problem (\ref{stiff1})--(\ref{stiff2}) with stepsize $h=100/n$,  and {\tt ode15s} with the default parameters.}
\label{stiff}
\centerline{\begin{tabular}{|r|r|r|c|c|}
\hline
$n$ & $k$ & $s$ & time & $e_y$ \\
\hline
  50 & 40 & 38 &  0.09  & 2.92e-11  \\ 
  75 & 32 & 30 &  0.12  & 1.53e-11  \\ 
 100 & 28 & 26 &   0.17 & 1.93e-12  \\ 
  125 & 25 & 23 &  0.21  & 6.28e-12  \\ 
 150 & 22 & 20 &   0.27 & 9.43e-12  \\ 
\hline
\multicolumn{3}{|c|}{\tt ode15s} & 0.68  & 3.76e-04\\
\hline
\end{tabular}}
\end{table}

\subsubsection*{A stiff ODE-IVP}
At last, we consider a stiff ODE-IVP, 
\begin{equation}\label{stiff1}
\dot y(t) = \pmatrix{rrr}
       -9999 &          1    &       1\\
        9900  &      -100  &         1\\
          98    &      98     &     -2\endpmatrix \big[ y(t) - g(t)\big] + \dot{g}(t), \qquad y(0) = g(0),
\end{equation}
with $g(t)$ a known function, having evidently solution $y(t)=g(t)$. We choose
\begin{equation}\label{stiff2}
g(t) = \pmatrix{ccc} \cos 2\pi t, & \cos 4\pi t, &\cos 6\pi t\endpmatrix^\top,
\end{equation}
and consider the SHBVM with $tol\approx 10^{-8}$ in (\ref{newtol}) for choosing $s$ (as before, $\kappa$ is chosen according to (\ref{kappas})), so that full accuracy is expected  in the numerical solution. The time-step used is $h=100/n$ for $n$ steps. The measured errors in the last point (coinciding with the initial condition), are then reported in Table~\ref{stiff}. For comparison, also the results obtained by the Matlab solver {\tt ode15s}, using its default parameters, are listed in the table. This latter solver requires 6006 steps, with time-steps  approximately in the range  $[1.9\cdot10^{-3},\,2\cdot10^{-2}]$.

\section{Conclusions}\label{fine}
In this paper we provide a thorough analysis of SHBVMs, namely HBVMs used as spectral methods in time, which further confirms their effectiveness. From the analysis, one obtains that the super-convergence of HBVMs is maintained also when using relatively large time-steps. SHBVMs become a practical method, due to the very efficient nonlinear {\em blended} iteration inherited from HBVMs. As a consequence, SHBVMs appear to be good candidates as {\em general ODE solvers}. This is indeed confirmed by a few numerical tests concerning a Hamiltonian problem, a Poisson (not Hamiltonian) problem, and a stiff ODE-IVP. The same tests show the numerical assessment of the theoretical achievements.

\end{document}